\documentclass{article}
\usepackage{graphicx}
\title{Quantum Invariants of Links and New Quantum Field Models}

\author{Sze Kui Ng
\\  Department of Mathematics,
Hong Kong Baptist University, Hong Kong
\\E-mail: skng@hkbu.edu.hk
\\Tel: (852)2339 7018
\\Fax: (852)2339 5811
}
\begin{document}
\date{}
\maketitle
\begin{abstract}

We propose a gauge model of quantum
electrodynamics (QED) and its nonabelian generalization 
from which we derive 
knot invariants such as the Jones polynomial. Our approach is inspired by the work of Witten who derived knot invariants
from quantum field theory based on the Chern-Simon Lagrangian.
From our approach we can derive new knot and link invariants which extend
the Jones polynomial and give a complete classification of knots
and links. From these new knot invariants we have that knots can be
completely classified by the power index $m$ of $TrR^{-m}$
where $R$ denotes the $R$-matrix for braiding and is the monodromy 
of the Knizhnik-Zamolodchikov equation. A classification table of knots
can then be formed where prime knots are classified by prime integer $m$
and nonprime knots are classified by nonprime integer $m$.



\end{abstract}

\section{Introduction}\label{sec00}

In 1989 Witten derived knot invariants such as
the Jones polynomial from quantum field theory based on the
Chern-Simon Lagrangian \cite{Witten}.
Inspired by Witten's work
in this paper we shall derive knot invariants from a 
gauge model of Quantum electrodynamics (QED) and its nonabelian
generalization.
 From our approach we shall first derive the Jones
polynomial and then we derive new knot and link invariants which
extend the Jones polynomial and can give a complete classification of knots
and links. 

This paper is organized as follows. In section \ref{sec2}
we give a brief description of a gauge model
of QED and its nonabelian extention. In this paper we shall consider a 
nonabelian extension
with a $SU(2)$ gauge symmetry. With this quantum field model we
introduce the partition function and the correlation of a Wilson loop which will be a
knot invariant of the trivial knot (also called the unknot). This correlation of
Wilson loop will later be generalized to be knot invariants
of nontrivial knots and links. To investigate the properties of these
partition function and correlations
in section \ref{sec5} we derive a chiral symmetry from
the gauge transformation of this new quantum field model.
From this chiral symmetry in section \ref{sec6},
section \ref{sec7} we derive a conformal field theory
which contains topics such as the affine Kac-Moody algebra
and the Knizhnik-Zamolodchikov equation.
This KZ equation is an equation of correlations from
which in section \ref{sec9} we can derive the skein relation
of the Jones polynomial. A main point of our theory on the KZ equation is that we can derive two KZ equations
which are dual to each other. From these two KZ equations
we derive a quantum group structure for the $W$ matrices
from which a Wilson loop is formed. Then from the
correlation of these $W$ matrices in section \ref{sec10} and section \ref{sec11}
we derive new knot invariants which extend the Jones
polynomial and gives a complete classification of knots.
In section \ref{sec12} we extend the new invariants to the case of links and we compute these new link invariants for some
examples of links. Then in section \ref{sec13} with the new knot invariants we give a classification table
of knots  which is formed
by using the power index $m$ which comes from these new knot invariants of the form 
$Tr R^{-m}\langle W(z,z)\rangle$ where $W(z,z)$ denotes
a Wilson loop  and $R$ is the 
braiding
matrix for the quantum group structure and is the monodromy of the two KZ equations.

\section{New gauge Model of QED and Nonabelian Extensions}\label{sec2}

To begin our derivation of knot and link invariants let us 
first describe a quantum field model.
Similar to the Wiener measure for the Brownian motion
which is constructed from the
integral $\int_{t_0}^{t_1}\left(\frac{dx}{dt}\right)^2dt$
we construct a measure for QED from the 
following energy integral:
\begin{equation}
-\frac12\int_{s_0}^{s_1}[
\frac12\left(\frac{dA_1}{ds}-\frac{dA_2}{ds}\right)^2
+\sum_{i=1}^2
\left(\frac{dz}{ds}-ieA_iz\right)^*\left(\frac{dz}{ds}-ieA_iz\right)]ds
\label{1.1}
\end{equation}
where $s$ denotes the proper time in relativity; $e$
denotes the electric charge and the complex variable $z$, real
variables $A_1$, $A_2$ represent one electron and two photons
respectively. 
By extending $ds$ to $(ih+\beta)ds$ we get a quantum theory of QED
where $h>0$ denotes the Planck constant and $\beta>0$ is a constant
related to absolute temperature.

The integral (\ref{1.1}) has the following gauge symmetry:
\begin{equation}
z'(s) = z(s)e^{iea(s)}, \quad
A'_i(s) = A_i(s)+\frac{da}{ds}
\quad i=1,2
\label{1.2}
\end{equation}
where $a(s)$ is a real valued function.

We remark that a main feature of (\ref{1.1}) is that
it is not formulated with the four-dimensional space-time but is formulated with the one dimensional proper time.
We refer to \cite{Ng} for the physical motivation of this new
QED theory.  

We can generalize the above QED model with $U(1)$ gauge symmetry
to QCD type models with nonabelian gauge symmetry.
As an illustration let us consider $SU(2)$ gauge symmetry.
Similar to 
(\ref{1.1}) we consider the following energy integral:
\begin{equation}
L := -\frac12\int_{s_0}^{s_1}
[\frac12 tr (D_1A_2-D_2A_1)^{*}(D_1A_2-D_2A_1) +
(D_1Z)^{*}(D_1Z)+(D_2Z)^{*}(D_2Z)]ds
\label{n1}
\end{equation}
where 
$Z= (z_1, z_2)^{T}$ is a two dimensional complex vector;
$A_j =\sum_{k=1}^{3}A_j^k t^k $ $(j=1,2)$ where
$A_j^k$ denotes a real component of a gauge field $A^k$;
$t^k$ denotes a generator of $SU(2)$;
and
$D_j=\frac{d}{ds}-igA_j$ $(j=1,2)$
where $g$ denotes the charge of interaction.
From (\ref{n1}) we can develop a QCD type model as similar
to that for the QED model.
We have that (\ref{n1}) is invariant under the following
gauge transformation:
\begin{equation}
\begin{array}{rl}
Z^{\prime}(s) &=U(a(s))Z(s) \\
A_j^{\prime}(s) &=U(a(s))A_j(s)U^{-1}(a(s))+U(a(s))\frac{dU^{-1}(a(s))}{ds},
j =1,2
\end{array}
\label{n2}
\end{equation}
where $U(a(s))=e^{-a(s)}$ and $a(s)=\sum_k a^k (s)t^k$.
We shall mainly consider the case that $a$ is a function
of the form $a(s)=\omega_1(r(s))$ where $\omega_1$ 
and $r$ are analytic functions.

\section{Knot Invariants} \label{sec4}

Since (\ref{n1}) is not formulated with the space-time, 
as analogous to the approach of Witten on knot invariants
\cite{Witten} the
following partition function will be shown to be a topological invariant for
knots:
\begin{equation}
\langle TrW_{R} \rangle := \int DA_1DA_2DZ
  e^{L} TrW_{R}(C)
\label{n3}
\end{equation}
where 
\begin{equation}
W_R(C):= W(r_0, r_1):= Pe^{\int_C A_idx^i}
\label{n4}
\end{equation}
which may be called a Wilson loop as analogous to the usual
Wilson loop 
\cite{Witten} where $C$ denotes a closed curve of the following
form
\begin{equation}
C(s) =(x^1(r(s)), x^2(r(s))), s_0\leq s \leq s_1
\label{n4a}
\end{equation}
where $r$ is an analytic function
such that $r_0 :=r(s_0)=r(s_1) :=r_1$. This closed curve
$C$ is
in a two dimensional phase plane $(x^1, x^2)$ which is 
dual to $(A_1, A_2)$.
We let this closed curve $C$ represents the
projection of a knot in this two dimensional space.
As usual the notation
$P$ in the definition of $W_R(C)$ denotes a path-ordered product and $R$ denotes a representation
of $SU(2)$ \cite{Kau}\cite{Baez}.

We remark that we also extend the definition of
$W_R(C)$ to the case that $C$ is not a closed curve
with $r_0\neq r_1$.

We shall show that (\ref{n3}) is a topological invariant
for a trivial knot. This means that in (\ref{n3}) the
closed curve $C$ represents a trivial knot.
We shall extend (\ref{n3}) to let $C$ represent nontrivial
knot.

Our aim is to compute the above knot invariant 
and its generalization to knot invariants of nontrivial knots
which will be defined.

\section{Chiral Symmetry} \label{sec5}

For a given curve $C(s)=(x^1(r(s)), x^2(r(s))), s_0 \leq s\leq s_1$ which may not be a closed curve
we define $W(r_0, r_1)$ by (\ref{n4}) where $r_0=r(s_0)$
and $r_1=r(s_1)$.
Then under an analytic gauge transformation
we have the following chiral symmetry:
\begin{equation}
W(r_0, r_1) \mapsto U(\omega(r_1))
W(r_0, r_1)U^{-1}(\omega(r_0))
\label{n5a}
\end{equation}
where $\omega$ denotes an analytic function.
This chiral symmetry is analogous to the chiral symmetry
of the usual nonabelian guage theory where $U$ denotes an element of $SU(2)$ \cite{Kau}.
We may extend (\ref{n5}) by extending $r$ to complex
variable $z$ to have the following chiral symmetry:
\begin{equation}
W(z_0, z_1) \mapsto U(\omega(z_1))W(z_0, z_1)
U^{-1}(\omega(z_0))
\label{n5}
\end{equation}
This analytic continuation corresponds to the complex
transformation $s \mapsto (ih+\beta)s$ for describing
quantum physics.

From this chiral symmetry  we have the following formulas for the
variations $\delta_{\omega}W$ and $\delta_{\omega^{\prime}}W$ with
respect to the chiral symmetry:
\begin{equation}
\delta_{\omega}W(z,z')=W(z,z')\omega(z)
\label{k1}
\end{equation}
and
\begin{equation}
\delta_{\omega^{\prime}}W(z,z')=-\omega^{\prime}(z')W(z,z')
\label{k2}
\end{equation}
where $z$ and $z'$ are independent variables and
$\omega^{\prime}(z')=\omega(z)$ when
$z'=z$. In (\ref{k1}) the variation is with respect to the
$z$ variable while in (\ref{k2}) the variation is with
respect to the $z'$ variable. This two-side-variations is possible
when $z\neq z'$.

\section{Affine Kac-Moody Algebra} \label{sec6}

Let us define
\begin{equation}
J(z) := -k W^{-1}(z, z')\partial_z W(z, z')
\label{n6}
\end{equation}
where $ k>0 $ is a constant.
As analogous to the WZW model \cite{Kni}\cite{Fra}
 $J$ is a generator of the chiral symmetry for (\ref{k1}).

Let us consider the following correlation
\begin{equation}
\langle W_{R}A(z) \rangle := \int DA_1DA_2DZ
  e^{L} W_{R}(C)A(z)
\label{n8a}
\end{equation}
By taking a gauge transformation on this correlation
and by the gauge invariance of (\ref{n1})
we can derive a Ward identity from which we have the following
relation:
\begin{equation}
\delta_{\omega}A(z)=\frac{-1}{2\pi i}\oint_z dw\omega(w)J(w)A(z)
\label{n8b}
\end{equation}
where $\delta_{\omega}A $ denotes the variation of the field $A$ with respect
to the chiral symmetry and the closed line integral $\oint$ is with
center $z$ and we let the generator $J$ be given by (\ref{n6}). We remark that our approach here is
analogous to the WZW model in conformal field theory
\cite{Kni}.

From (\ref{n5}) and (\ref{n6}) we have that the variation $\delta_{\omega}J$
of the generator $J$ of the chiral symmetry is given by
\cite{Kni}\cite{Fra}:
\begin{equation}
\delta_{\omega}J= \lbrack J, \omega\rbrack -k\partial_z \omega
\label{n8c}
\end{equation}

From (\ref{n8b}) and (\ref{n8c}) we have that $J$ satisfies
the following relation of current algebra \cite{Kni}\cite{Fra}\cite{Fuc}:
\begin{equation}
J^a(w)J^b(z)=\frac{k\delta_{ab}}{(w-z)^2}
+\sum_{c}if_{abc}\frac{J^c(z)}{(w-z)}
\label{n8d}
\end{equation}
where we write
\begin{equation}
J(z) = \sum_a J^a(z) t^a = \sum_a
\sum_{n=-\infty}^{\infty}J_n^a z^{-n-1} t^a
\label{n7}
\end{equation}
Then from (\ref{n8d}) we can show that $J_n^a$ satisfy the following affine
Kac-Moody algebra \cite{Kni}\cite{Fra}\cite{Fuc}:
\begin{equation}
[J_m^a, J_n^b] =
if_{abc}J_{m+n}^c + km\delta_{ab}\delta_{m+n, 0}
\label{n8}
\end{equation}
where the constant $k$ is called the central extension
or the leval of the Kac-Moody algebra.

Let us consider another generator of the chiral symmetry for (\ref{k2}) given by
\begin{equation}
J^{\prime}(z')= k\partial_{z'}W(z, z')W^{-1}(z, z')
\label{d1}
\end{equation}
Similar to $J$ by the following correlation:
\begin{equation}
\langle A(z')W_{R} \rangle := \int DA_1DA_2DZ
  A(z')W_{R}(C)e^{L}
\label{n8aa}
\end{equation}
we have the following formula for $J^{\prime}$:
\begin{equation}
\delta_{\omega^{\prime}}A(z')=
\frac{-1}{2\pi i}\oint_{z^{\prime}}^{-} dwA(z')J^{\prime}(w)
(-\omega^{\prime})(w)
=\frac{-1}{2\pi i}\oint_{z^{\prime}} dwA(z')J^{\prime}(w)
\omega^{\prime}(w)
\label{n8b1}
\end{equation}
where $\oint^{-}$ denotes an integral with clockwise
direction while $\oint$ denotes an integral with counterclockwise direction. We remark that this two-side variation from (\ref{n8a}) and (\ref{n8aa}) is important
for deriving the two KZ equations which are dual to each
other.

Then similar to (\ref{n8c}) we also we have
\begin{equation}
\delta_{\omega^{\prime}}J^{\prime}= 
\lbrack \omega^{\prime}, J^{\prime}\rbrack -k\partial_{z'} \omega^{\prime}
\label{n8c1}
\end{equation}
Then from (\ref{n8b1}) and (\ref{n8c1}) we can derive the current
algebra and the Kac-Moody algebra for $J^{\prime}$ which are of the
same form of (\ref{n8d}) and (\ref{n8}).

\section{Dual Knizhnik-Zamolodchikov Equation} \label{sec7}

Let us first consider (\ref{k1}).
From (\ref{n8b}) and (\ref{k1}) we have
\begin{equation}
J^a(z)W(w, w') \sim \frac{-t^aW(w,w')}{z-w}
\label{k3}
\end{equation}

Let us define an energy-momentum tensor $T(z)$ by
\begin{equation}
T(z) := \frac{1}{k+g}\sum_a :J^a(z)J^a(z):
\label{k4}
\end{equation}
where $g$ is the dual Coxter number. In (\ref{k4})
the symbol $:...:$ denotes normal ordering. This is the
Sugawara construction of energy-momentum tensor
where the appearing of $g$ is from a renormalization
of quantum effect by requiring the operator product
expansion of $T$ with itself to be of the following
form \cite{Kni} \cite{Fra} \cite{Fuc}:
\begin{equation}
T(z)T(w)=\frac{c}{2(z-w)^4}+
         \frac{2T(w)}{(z-w)^2}+\frac{\partial T(w)}{(z-w)}
\label{k5}
\end{equation}
for some constant $c=\frac{kd}{k+g}$ where $d$ denotes the dimension of $SU(2)$.

Then we have the following $TW$ operator product:
\begin{equation}
T(z)W(w,w')\sim \frac{\Delta}{(z-w)^2}
                +\frac{1}{(z-w)}L_{-1}W(w,w')
\label{k6}
\end{equation}
where $L_{-1}W(w,w')=\partial_{w}W(w,w')$ and
\begin{equation}
\Delta=\frac{\sum_a t^a t^a}{2(k+g)}=\frac{N^2-1}{2N(k+N)}
\label{k6a}
\end{equation}
where $N$ is for $SU(N)$.

From (\ref{k4}) and (\ref{k6}) we have the following equation \cite{Fra}\cite{Fuc}:
\begin{equation}
L_{-1}W(w,w')=\frac{1}{k+g}J_{-1}^aJ_{0}^aW(w,w')
\label{k7}
\end{equation}
Then form (\ref{k3}) we have 
\begin{equation}
J_{0}^aW(w,w')=-t^aW(w,w')
\label{k8}
\end{equation}
By (\ref{k7}) and (\ref{k8}) we have
\begin{equation}
\partial_z W(z, z')=\frac{-1}{k+g}J_{-1}^at^aW(z,z')
\label{k9}
\end{equation}

From this equation and by the $JW$ operator product
(\ref{k3}) we have the following
Knizhnik-Zamolodchikov equation \cite{Fra} \cite{Fuc}:
\begin{equation}
\partial_{z_i}
\langle W(z_1, z_1^{\prime})\cdot\cdot\cdot 
W(z_n, z_n^{\prime})
\rangle
=-\frac{1}{k+g}
\sum_{j\neq i}^{n}\frac{\sum_a t^a \otimes t^a}{z_i-z_j}
\langle W(z_1, z_1^{\prime})\cdot\cdot\cdot 
W(z_n, z_n^{\prime})
\rangle
\label{n9}
\end{equation}
We remark that in deriving (\ref{n9}) we have used
line integral expression of operators with counterclockwise
direction \cite{Fra}\cite{Fuc}.

 It is interesting and important that we also have
another Knizhnik-Zamolodchikov equation
which will be called the dual equation of
(\ref{n9}). The derivation of this dual equation is dual to the
above derivation in that the line integral for this derivation of dual
equation
is with clockwise direction in contrast to counterclockwise
direction in the above derivation and that the operator products and their corresponing variables are with reverse
order to that in the above derivation. 

From (\ref{k2}) and (\ref{n8b1}) we have a $WJ^{\prime}$ operator product given
by
\begin{equation}
W(w, w')J^{\prime a}(z') \sim \frac{-W(w, w')t^a}{w'-z'}
\label{d2}
\end{equation}
Similar to the above derivation of the KZ equation from (\ref{d2})
we can then derive the following Knizhnik-Zamolodchikov
equation which is dual to (\ref{n9}):
\begin{equation}
\partial_{z_i^{\prime}}
\langle W(z_1,z_1^{\prime})\cdot\cdot\cdot W(z_n,z_n^{\prime})
\rangle
= -\frac{1}{k+g}\sum_{j\neq i}^{n}
\langle W(z_1, z_1^{\prime})\cdot\cdot\cdot 
W(z_n, z_n^{\prime})
\rangle
\frac{\sum_a t^a\otimes t^a}{z_j^{\prime}-z_i^{\prime}}
\label{d8}
\end{equation}

\section{Skein Relation for Jones Polynomial}\label{sec9}

Following the idea of Witten \cite{Witten}, if we cut a knot we
get two pieces of curves crossing (or not crossing) each other
once. This gives two primary fields
$W(z_1, z_2)$ and $W(z_3, z_4)$ where $W(z_1, z_2)$ corresponds
to a piece of curve with end points parametrized by $z_1$ and
$z_2$ and $W(z_3, z_4)$ corresponds to the other piece of curve
with end points parametrized by $z_3$ and $z_4$.
Let us write
\begin{equation}
W(z_i, z_j) = W(z_i, z_i^{\prime})W^{-1}(z_j, z_j^{\prime})
\label{n10}
\end{equation}
for $i=1,3$ and $j=2, 4$ and for some $z_k^{\prime}$
with $z_1^{\prime}=z_2^{\prime}$ and
$z_3^{\prime}=z_4^{\prime}$.
These two pieces of curves then correspond to the following
four-point correlation function:
\begin{equation}
G(z_1, z_2, z_3, z_4) :=
\langle W(z_1^{\prime}, z_1)W^{-1}(z_2^{\prime}, z_2)
W(z_3, z_3^{\prime})W^{-1}(z_4, z_4^{\prime})\rangle
\label{n11}
\end{equation}
(In the notation $G(z_1,z_2,z_3,z_4)$ we have supressed 
the $z'$ variables for simplicity). Then we
have \cite{Fra}\cite{Fuc}:
\begin{equation}
G(z_1,z_2,z_3,z_4)=
[(z_1-z_3)(z_2-z_4)]^{-2\Delta}G(x)
\label{n12}
\end{equation}
where $\Delta =\frac{N^2-1}{2N(N+k)}$ and
$x= \frac{(z_1-z_2)(z_3-z_4)}{(z_1-z_3)(z_2-z_4)}$ and
from the KZ equation $G(x)$ satisfies the following
equation:
\begin{equation}
\frac{dG}{dx}= [\frac{1}{x}P +\frac{1}{x-1} Q]G
\label{n13}
\end{equation}
where
\begin{equation}
P= -\frac{1}{N(N+k)}\left(\begin{array}{cc}
  N^2-1 & N \\
  0     & -1 \end{array}\right),
\quad
Q= -\frac{1}{N(N+k)}\left(\begin{array}{cc}
   -1 & 0 \\
   N  & N^2-1\end{array}\right)
\label{n14}
\end{equation}
This equation has two independent conformal block solutions
forming a vector space of dimension 2.
Let $\psi$ be a vector in this space and let $B$ denotes the
braid operation. Then following Witten
\cite{Witten} we have
\begin{equation}
 a\psi -bB\psi +B^2\psi=0
\label{n15}
\end{equation}
where $a=det B$ and $b= Tr B$. Then following Witten
\cite{Witten} from (\ref{n15}) we can derive the following skein
relation for the Jones polynomial:
\begin{equation}
\frac1{t}V_{L_-} -tV_{L_+}
= (t^{\frac12} - \frac{1}{t^{\frac12}})V_{L_0}
\label{n16}
\end{equation}
where $V_{L_-}$, $V_{L_+}$ and $V_{L_0}$ are the Jones polynomials
for undercrossing, overcrossing and zero crossing respectively.

\section{New Knot Invariants Extending Jones Polynomial}\label{sec10}

Let us consider again the correlation
$G(z_1, z_2, z_3, z_4)$ in (\ref{n11}) which also have the
following form:
\begin{equation}
G(z_1, z_2, z_3, z_4)=
\langle W(z_1, z_2)W(z_3, z_4)\rangle
\label{m1}
\end{equation}
From it in this section we shall present a method which is
different from the above section to derive new knot invariants. These new knot invariants will
extend the Jones polynomial and they will be defined by
generalizing (\ref{n3}).

We have that $G$ satisfies the KZ equation for the
variables $z_1$, $z_3$ and satisfies the dual KZ equation
for the variables $z_2$ and $z_4$.
By solving the KZ equation we have that $G$ is of the
form 
\begin{equation}
e^{t\log (z_1-z_3)}C_1
\label{m2}
\end{equation}
where $t:=\frac{1}{k+g}\sum_a t^a \otimes t^a$
and $C_1$ denotes a constant matrix which is independent
of the variable $z_1-z_3$.

Similarly by solving the dual KZ equation we have that
$G$ is of the form
\begin{equation}
C_2e^{-t\log (z_4-z_2)}
\label{m3}
\end{equation}
where $C_2$ denotes a constant matrix which is independent
of the variable $z_4-z_2$.

From (\ref{m2}), (\ref{m3}) and we let
$C_1=Ae^{-t\log(z_4-z_2)}$, $C_2= e^{t\log(z_1-z_3)}A$ where $A$ is a constant matrix we have that
$G$ is given by
\begin{equation}
G(z_1, z_2, z_3, z_4)=
e^{t\log (z_1-z_3)}Ae^{-t\log (z_4-z_2)}
\label{m4}
\end{equation}

Now let $z_2=z_3$.
Then as $z_4 \to z_1$ we have
\begin{equation}
TrG(z_1, z_2, z_2, z_1)=
Tre^{i2n\pi t}A
=:TrR^{2n}A \quad\quad n=0, \pm 1, \pm 2, ...
\label{m5}
\end{equation}
where $R=e^{i\pi t}$ is the monodromy of the the KZ equation \cite{Chari}. We remark that (\ref{m5}) is a
multivalued function.
 From (\ref{m5}) we have the following relation between
the partition function $Z$ and the matrix $A$:
\begin{equation}
 A=IZ
\label{m5a}
\end{equation}
where $I$ denotes the identity matrix.
Now let $C$ be a closed curve in the complex plane
with initial and final end points $z_1$. Then the
following correlation function
\begin{equation}
Tr\langle W(z_1, z_1)\rangle
=Tr\langle W(z_1, z_2)W(z_2, z_1)\rangle
\label{m6}
\end{equation}
which is the definition (\ref{n3}) defined along the
curve $C$, with $W(z_1, z_1)=W(z_1,z_2)W(z_2,z_1)$,
can be regarded as a knot invariant of the trivial
knot in the three dimensional space whose porjection
in the complex plane is the curve $C$.
Indeed, from (\ref{m1}) and (\ref{m5}) we can compute
(\ref{m6}) which is given by:
\begin{equation}
Tr\langle W(z_1,z_1)\rangle 
=ZTr R^{2n} \quad\quad n=0,\pm 1, \pm 2, ...
\label{m6a}
\end{equation}
From (\ref{m6a}) we see that (\ref{m6}) is
independent of the closed curve $C$ which represents
the projection of a trivial knot and thus can be
regarded as a knot invariant for the trivial knot.

In the following let us extend the definition (\ref{m6})
to knot invariants for nontrivial knots.

Since $R $ is the monodromy
of the KZ equation, we have a branch cut such that
\begin{equation}
\langle W(z_3,w)W(w, z_2)W(z_1,w)W(w,z_4)\rangle =
R\langle W(z_1,w)W(w,z_2)W(z_3,w)W(w,z_4)\rangle
\label{m7}
\end{equation}
where $z_1$ and $z_3$ denote two points on a closed curve
such that along the direction of the curve the point
$z_1$ is before the point $z_3$. From (\ref{m7}) we have
\begin{equation}
W(z_3,w)W(w,z_2)W(z_1,w)W(w,z_4)
=RW(z_1,w)W(w,z_2)W(z_2,w)W(w,z_4)
\label{m7a}
\end{equation}
Similarly for the dual KZ equation we have
\begin{equation}
\langle W(z_1,w)W(w, z_4)W(z_3,w)W(w,z_2)\rangle =
\langle W(z_1,w)W(w,z_2)W(z_3,w)W(w,z_4)\rangle R^{-1}
\label{m8}
\end{equation}
and 
\begin{equation}
W(z_1,w)W(w,z_4)W(z_3,w)W(w,z_2)=
W(z_1,w)W(w,z_2)W(z_3,w)W(w,z_4)R^{-1}
\label{m8a}
\end{equation}
where $z_2$ before $z_4$.
From (\ref{m7a}) and (\ref{m8a}) we have
\begin{equation}
 W(z_3,z_4)W(z_1,z_2)=
RW(z_1,z_2)W(z_3,z_4)R^{-1}
\label{m9}
\end{equation}
where $z_1$ and $z_2$ denote the end points of a curve which is before a curve with end points $z_3$ and $z_4$.
From (\ref{m9}) we see that the algebraic structure of these $W$
matrices is analogous to the quasi-triangular quantum
group \cite{Chari}\cite{Fuc}.
Now we let $W(z_i,z_j)$ represent a piece of curve
with initial end point $z_i$ and final end point $z_j$.
Then we let 
\begin{equation}
W(z_1,z_2)W(z_3,z_4)
\label{m11}
\end{equation}
represent two pieces of uncrossing curve.
Then by interchanging $z_1$ and $z_3$ we have
\begin{equation}
W(z_3,w)W(w,z_2)W(z_1,w)W(w,z_4)
\label{m12}
\end{equation}
represent the curve specified by $W(z_1,z_2)$ upcrossing the
curve specified by $W(z_3,z_4)$ at $z$.
Similarly by interchanging $z_2$ and $z_4$ we have
\begin{equation}
W(z_1,w)W(w,z_4)W(z_3,w)W(w,z_2)
\label{m13}
\end{equation}
represent the curve specified by $W(z_1,z_2)$ undercrossing
the curve specified by $W(z_3,z_4)$ at $z$.

Now for a closed curve we may cut it into a sum of
parts which are formed by two pieces of curve crossing or
not crossing each other. Each of these parts is represented
by (\ref{m11}), (\ref{m12}) or (\ref{m13}).
Then we may define a correlation for a knot whose projection
is this closed curve by the following form:
\begin{equation}
 Tr\langle \cdot\cdot\cdot 
 W(z_3,z)W(z,z_2)W(z_1,z)W(z,z_4)
\cdot\cdot\cdot 
 \rangle
\label{m14}
\end{equation}
where we use (\ref{m12}) as an example to represent the state of the
two pieces of curve specified by 
 $W(z_1,z_2)$ and
$W(z_3,z_4)$. The 
 $\cdot\cdot\cdot$ means multiplications
of a sequence of parts represented by (\ref{m11}),
(\ref{m12}) or (\ref{m13}) according to the state of
each part. The order of the sequence in (\ref{m14})
 follows the order of the parts given by the direction of the
knot.

We shall show that 
(\ref{m14}) is  a knot invariant for a given knot. In the following let us consider
some examples to illustrate the way to define (\ref{m14}) and to show that (\ref{m14}) is a knot invariant. We shall also derive the three Reidemeister moves for
the equivalence of knots.

Let us first consider a knot in Fig. 1.
For this knot we have that (\ref{m14}) is given by
\begin{equation}
Tr\langle W(z_2,w)W(w,z_2)W(z_1,w)W(w,z_1)\rangle
\label{m15a}
\end{equation}
where the product of $W$  is from the definition (\ref{m12}).
In applying (\ref{m12}) we let $z_1$ as the
starting and the final point. We remark that the
$W$ matrices  must be put together
to follow the definition (\ref{m12}) and they are
not separated to follow the direction of the knot.

\begin{figure}[hbt]
\centering
\includegraphics[scale=0.6]{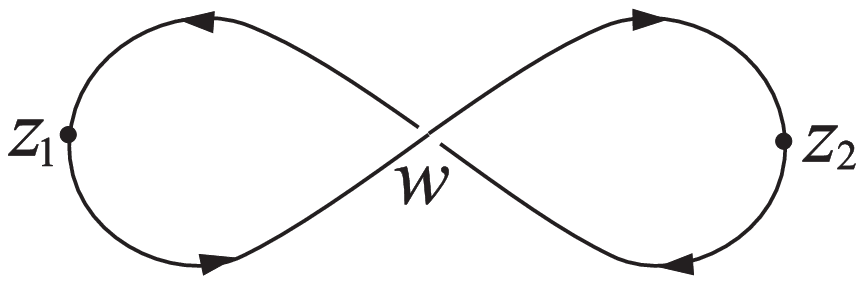}
                  
Fig.1
\end{figure}
Then  we have that (\ref{m15a}) is equal to
\begin{equation}
\begin{array}{rl}
&Tr\langle W(w,z_2)W(z_1,w)W(w,z_1)W(z_2,w)\rangle \\
=&Tr\langle RW(z_1,w)W(w,z_2)R^{-1}
RW(z_2,w)W(w,z_1)R^{-1}\rangle \\
=&Tr\langle W(z_1,z_2)W(z_2,z_1)\rangle \\
=&Tr\langle W(z_1,z_1)\rangle 
\end{array}
\label{m16}
\end{equation}
where we have used (\ref{m9}).
We see that (\ref{m16}) is just the knot invariant (\ref{m6}) of a trivial knot.
Thus the knot in Fig.1 is with the same knot invariant of a trivial knot and this agrees with the fact that this knot is topologically equivalent
to a trivial knot.

Then let us derive the Reidemeister move 1. Consider
the diagram in Fig.2. We have that by (\ref{m12}) the definition (\ref{m14})
 for this
diagram is given by:
\begin{equation}
\begin{array}{rl}
& Tr\langle W(z_2,w)W(w,z_2)W(z_1,w)W(w,z_3)\rangle\\
=& Tr\langle W(z_2,w)RW(z_1,w)W(w,z_2)
R^{-1}W(w,z_3)\rangle\\
=& Tr\langle W(z_2,w)RW(z_1,z_2)R^{-1}W(w,z_3)\rangle\\
=& Tr\langle R^{-1}W(w,z_3)W(z_2,w)RW(z_1,z_2)\rangle\\
=& Tr\langle W(z_2,w)W(w,z_3)W(z_1,z_2)\rangle\\
=& Tr\langle W(z_2,z_3)W(z_1,z_2)\rangle\\
=& Tr\langle W(z_1,z_3)\rangle\\
\end{array}
\label{m18}
\end{equation}
where$W(z_1,z_3)$ represent a piece of curve with initial
end point $z_1$ and final end point $z_3$ which has no
crossing. When Fig.2 is a part of a knot we can also
derive a result similar to (\ref{m18}) which is for the
Reidemeister move 1. This shows that the Reidemeister move 1 holds. 
\begin{figure}[hbt]
\centering
\includegraphics[scale=0.6]{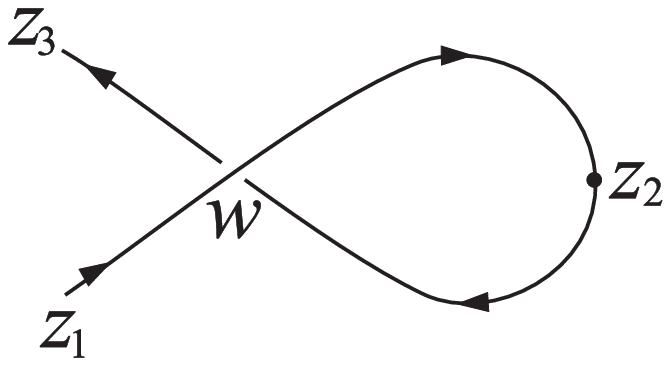}

                    Fig.2
\end{figure}

Then let us derive Reidemeister move 2.
By (\ref{m12}) we have that the definition (\ref{m14})
for the two pieces of curve in 
Fig.3a is given by
\begin{equation}
Tr\langle W(z_5,w_1)W(w_1,z_2)W(z_1,w_1)W(w_1,z_6)\cdot
W(z_4,w_2)W(w_2,z_3)W(z_2,w_2)W(w_2,z_5)\rangle
\label{m19}
\end{equation}
where the two products of $W$ separated by the 
 $\cdot$
are for the two crossings in Fig.3a.
We have that (\ref{m19}) is
equal to
\begin{equation}
\begin{array}{rl}
&Tr\langle W(z_4,w_2)W(w_2,z_3)W(z_2,w_2)W(w_2,z_5)\cdot
          W(z_5,w_1)W(w_1,z_2)W(z_1,w_1)W(w_1,z_6)\rangle \\
= &Tr\langle W(z_4,w_2)W(w_2,z_3)W(z_2,w_2)W(w_2,w_1)
             W(w_1,z_2)W(z_1,w_1)W(w_1,z_6)\rangle \\
=&Tr\langle W(z_4,w_2)W(w_2,z_3)W(z_2,w_2)RW(w_1,z_2)W(w_2,w_1)R^{-1}
            W(z_1,w_1)W(w_1,z_6)\rangle \\
=&Tr\langle W(z_4,w_2)RW(z_2,w_2)W(w_2,z_3)W(w_1,z_2)W(w_2,w_1)R^{-1}
            W(z_1,w_1)W(w_1,z_6)\rangle \\
=&Tr\langle W(z_4,w_2)RW(z_2,z_3)W(w_1,z_2)W(w_2,w_1)R^{-1}
            W(z_1,w_1)W(w_1,z_6)\rangle \\
=&Tr\langle W(z_4,w_2)W(w_1,z_2)W(z_2,z_3)RW(w_2,w_1)R^{-1}
            W(z_1,w_1)W(w_1,z_6)\rangle \\
=&Tr\langle W(z_4,w_2)W(w_1,z_3)RW(w_2,w_1)R^{-1}
            W(z_1,w_1)W(w_1,z_6)\rangle \\
=&Tr\langle RW(w_1,z_3)W(z_4,w_2)W(w_2,w_1)R^{-1}
            W(z_1,w_1)W(w_1,z_6)\rangle \\
=&Tr\langle RW(w_1,z_3)W(z_4,w_1)R^{-1}
            W(z_1,w_1)W(w_1,z_6)\rangle \\
=&Tr\langle W(z_4,w_1)W(w_1,z_3)
            W(z_1,w_1)W(w_1,z_6)\rangle \\
=&Tr\langle W(z_4,w_1)RW(z_1,w_1)W(w_1,z_3)
R^{-1}W(w_1,z_6)\rangle \\
=&Tr\langle W(z_1,z_3)R^{-1}W(w_1,z_6)W(z_4,w_1)R\rangle \\
=&Tr\langle W(z_1,z_3)W(z_4,w_1)W(w_1,z_6)\rangle \\
=&Tr\langle W(z_1,z_3)W(z_4,z_6)\rangle
\end{array}
\label{m20}
\end{equation}
where we have used (\ref{m9}). This shows that the diagram in Fig.3a is equivalent to
two uncrossing curves.
When Fig.3a is a part of a knot
we can also derive a result similar to (\ref{m20}) for
the Reidemeister
move 2.
This shows that the Reidemeister
move 2 holds. 
\begin{figure}[hbt]
\centering
\includegraphics[scale=0.5]{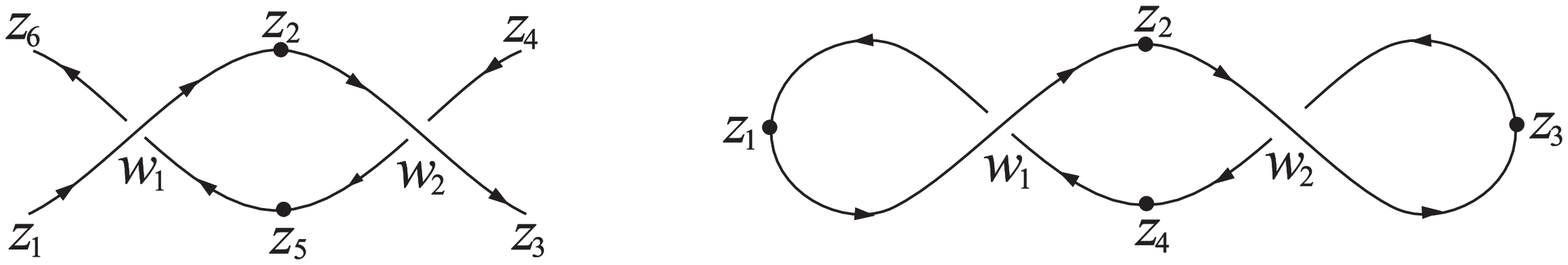}

               Fig.3a   \hspace*{5cm}          Fig.3b
\end{figure}
As an illustration
let us consider the knot in Fig. 3b which is related to
the Reidemeister move 2. By (\ref{m12}) we have that
the definition (\ref{m14}) for this knot is given by
\begin{equation}
\begin{array}{rl}
&Tr\langle W(z_3,w_2)W(w_2,z_3)W(z_2,w_2)W(w_2,z_4)\cdot
W(z_4,w_1)W(w_1,z_2)W(z_1,w_1)W(w_1,z_1)\rangle \\
=&Tr\langle RW(z_2,w_2)W(w_2,z_3)W(z_3,w_2)W(w_2,z_4)\cdot
W(z_4,w_1)W(w_1,z_1)W(z_1,w_1)W(w_1,z_2)R^{-1}\rangle \\
=&Tr\langle W(z_2,w_2)W(w_2,z_3)W(z_3,w_2)W(w_2,z_4)\cdot
W(z_4,w_1)W(w_1,z_1)W(z_1,w_1)W(w_1,z_2)\rangle \\
=&Tr\langle W(z_2,z_2)\rangle 
\end{array}
\label{m20a}
\end{equation}
where we let the curve be with $z_2$ as the initial
and final end point and we have used (\ref{m7a}) and
(\ref{m8a}). This shows that the knot in Fig.3b is with the same knot
invariant of a trivial knot.
This agrees with the fact
that this knot is equivalent to the trivial knot.

Let us then consider a trefoil knot in Fig.4a.
By (\ref{m12}) and similar to the above examples
we have that the definition (\ref{m14})
for this knot is given by:
\begin{equation}
\begin{array}{rl}
&Tr\langle W(z_4,w_1)W(w_1,z_2)W(z_1,w_1)W(w_1,z_5)\cdot
W(z_2,w_2)W(w_2,z_6) \\
&W(z_5,w_2)W(w_2,z_3)\cdot 
W(z_6,w_3)W(w_3,z_4)W(z_3,w_3)W(w_3,z_1)\rangle \\
=&Tr\langle W(z_4,w_1)RW(z_1,w_1)W(w_1,z_2)
R^{-1}W(w_1,z_5)\cdot
W(z_2,w_2)RW(z_5,w_2) \\
&W(w_2,z_6)R^{-1}W(w_2,z_3)\cdot 
W(z_6,w_3)RW(z_3,w_3)W(w_3,z_4)R^{-1}W(w_3,z_1)\rangle \\
=&Tr\langle 
W(z_4,w_1)RW(z_1,z_2)R^{-1}W(w_1,z_5)\cdot
W(z_2,w_2)RW(z_5,z_6)R^{-1}W(w_2,z_3)\cdot \\
&W(z_6,w_3)RW(z_3,z_4)R^{-1}W(w_3,z_1)\rangle \\
=&Tr\langle 
W(z_4,w_1)RW(z_1,z_2)
W(z_2,w_2)W(w_1,z_5)W(z_5,z_6)R^{-1}W(w_2,z_3)\cdot \\
&W(z_6,w_3)RW(z_3,z_4)R^{-1}W(w_3,z_1)\rangle \\
=&Tr\langle 
W(z_4,w_1)RW(z_1,w_2)W(w_1,z_6)R^{-1}W(w_2,z_3) \\&
W(z_6,w_3)RW(z_3,z_4)R^{-1}W(w_3,z_1)\rangle \\
=&Tr\langle 
W(z_4,w_1)W(w_1,z_6)W(z_1,w_2)W(w_2,z_3) \\&
W(z_6,w_3)RW(z_3,z_4)R^{-1}W(w_3,z_1)\rangle \\
=&Tr\langle 
W(z_4,z_6)W(z_1,z_3)
W(z_6,w_3)RW(z_3,z_4)R^{-1}W(w_3,z_1)\rangle \\
=&Tr\langle R^{-1}W(w_3,z_1)
W(z_4,z_6)W(z_1,z_3)
W(z_6,w_3)RW(z_3,z_4)\rangle \\
=&Tr\langle 
W(z_4,z_6)W(w_3,z_1)R^{-1}W(z_1,z_3)
W(z_6,w_3)RW(z_3,z_4)\rangle \\
=&Tr\langle 
RW(z_3,z_6)W(w_3,z_1)R^{-1}W(z_1,z_3)
W(z_6,w_3)\rangle \\
=&Tr\langle 
W(w_3,z_1)W(z_3,z_6)W(z_1,z_3)
W(z_6,w_3)\rangle \\
=&Tr\langle 
W(z_6,z_1)W(z_3,z_6)W(z_1,z_3)
\rangle 
\end{array}
\label{m21}
\end{equation}
where we have repeatly used (\ref{m9}). Then
 we have that (\ref{m21}) is equal to:
\begin{equation}
\begin{array}{rl}
&Tr\langle
W(z_6,w_3)W(w_3,z_1)W(z_3,w_3)W(w_3,z_6)W(z_1,z_3)\rangle 
\\
=&Tr\langle
RW(z_3,w_3)W(w_3,z_1)W(z_6,w_3)W(w_3,z_6)W(z_1,z_3)\rangle\\
=&Tr\langle
RW(z_3,w_3)RW(z_6,w_3)W(w_3,z_1)
R^{-1}W(w_3,z_6)W(z_1,z_3)\rangle\\
=&Tr\langle
W(z_3,w_3)RW(z_6,z_1)
R^{-1}W(w_3,z_6)W(z_1,z_3)R\rangle\\
=&Tr\langle
W(z_3,w_3)RW(z_6,z_3)W(w_3,z_6)\rangle\\
=&Tr\langle W(w_3,z_6)W(z_3,w_3)RW(z_6,z_3)\rangle\\
=&Tr\langle RW(z_3,w_3)W(w_3,z_6)W(z_6,z_3)\rangle\\
=&Tr\langle RW(z_3,z_3)\rangle
\end{array}
\label{m22}
\end{equation}
where we have used (\ref{m7a}) and (\ref{m9}).
We see that (\ref{m22}) is a knot invariant for the trefoil knot in Fig.4a.

\begin{figure}[hbt]
\centering
\includegraphics[scale=0.5]{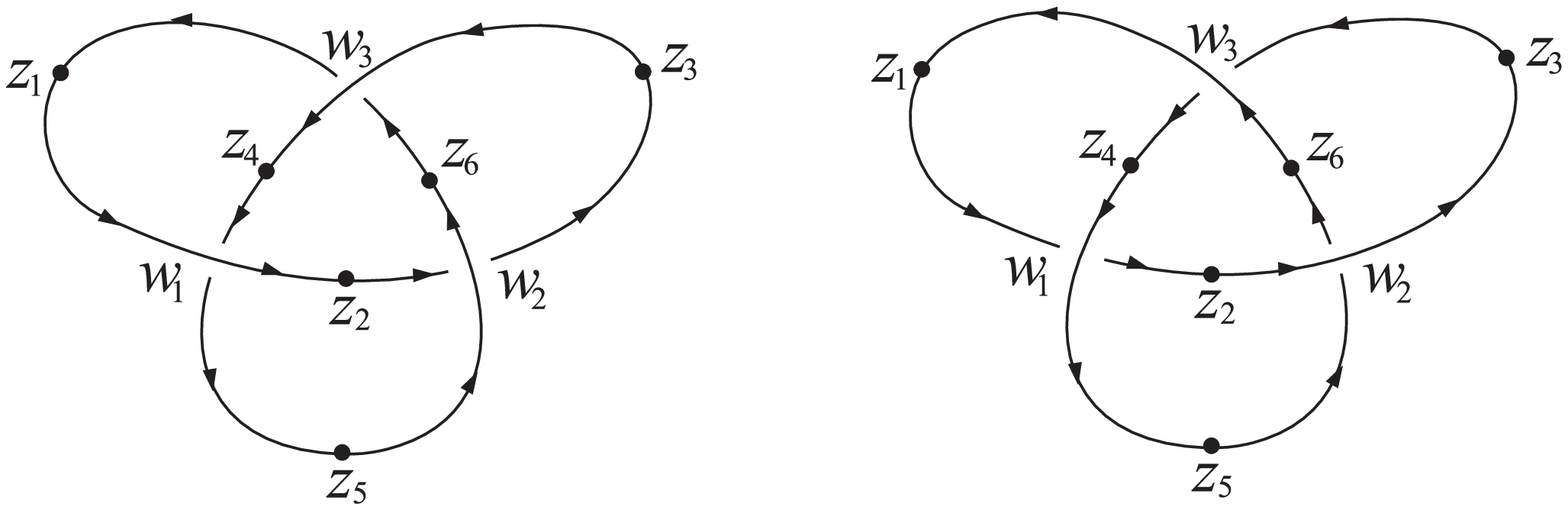}

             Fig.4a  \hspace*{5cm}  Fig.4b
\end{figure}

Then let us consider the trefoil knot in Fig. 4b which
is the mirror image of the trefoil knot in Fig.4a.
The definition (\ref{m14}) for this knot is given by:
\begin{equation}
\begin{array}{rl}
&Tr \langle W(z_1,w_1)W(w_1,z_5)W(z_4,w_1)W(w_1,z_2)\cdot
W(z_5,w_2)W(w_2,z_3)W(z_2,w_2)W(w_2,z_6)\cdot \\
&W(z_3,w_3)W(w_3,z_1)W(z_6,w_3)W(w_3,z_4)\rangle \\
=&Tr\langle W(z_5,z_1)W(z_2,z_5)W(z_1,z_2)\rangle
\end{array}
\label{m25}
\end{equation}
where similar to (\ref{m21}) we have repeatly used
(\ref{m9}).
Then we have that (\ref{m25}) is
equal to:
\begin{equation}
\begin{array}{rl}
&Tr\langle W(z_5,z_1)W(z_2,w_1)W(w_1,z_5)W(z_1,w_1)W(w_1,z_2)\rangle\\
=&Tr\langle 
W(z_5,z_1)W(z_2,w_1)W(w_1,z_2)
W(z_1,w_1)W(w_1,z_5)R^{-1}\rangle\\
=&Tr\langle 
W(z_5,z_1)W(z_2,w_1)RW(z_1,w_1)W(w_1,z_2)
R^{-1}W(w_1,z_5)R^{-1}\rangle\\
=&Tr\langle 
R^{-1}W(z_5,z_1)W(z_2,w_1)RW(z_1,z_2)
R^{-1}W(w_1,z_5)\rangle\\
=&Tr\langle 
W(z_2,w_1)W(z_5,z_2)
R^{-1}W(w_1,z_5)\rangle\\
=&Tr\langle 
W(z_5,z_2)
R^{-1}W(w_1,z_5)W(z_2,w_1)\rangle\\
=&Tr\langle 
W(z_5,z_2)
W(z_2,w_1)W(w_1,z_5)R^{-1}\rangle\\
=&Tr\langle 
W(z_5,z_5)R^{-1}\rangle
\end{array}
\label{m26}
\end{equation}
where we have used (\ref{m8a}) and (\ref{m9}).
We see that this is a knot invariant for the trefoil knot in Fig.4b. we notice that the knot invariants for the two
trefoil knots are different. This shows that these two
trefoil knots are not topologically equivalent.

Then let us derive the Reidemeister move 3. We have that
the definition (\ref{m14}) for the diagram in Fig.5a
is given by
\begin{equation}
\begin{array}{rl}
&Tr\langle W(z_7,w_1)W(w_1,z_2)W(z_1,w_1)W(w_1,z_8)\cdot
W(z_8,w_3)W(w_3,z_6)W(z_5,w_3)W(w_3,z_9)\cdot\\
&W(z_4,w_2)W(w_2,z_3)W(z_2,w_2)W(w_2,z_5)
\rangle
\end{array}
\label{m29}
\end{equation}
where we let the the curve with end points $z_7$, $z_9$
starts first, then the curve with end points $z_1$, $z_3$
starts second.
Similar to the derivation of the above invariants by (\ref{m7a}),
(\ref{m8a}), (\ref{m9}) we have that
(\ref{m29}) is equivalent to the following correlation:
\begin{equation}
Tr\langle W(z_7,z_9)W(z_1,z_3)W(z_5,z_6)
R^{-1}W(z_4,z_5)\rangle
\label{m30}
\end{equation}

\begin{figure}[hbt]
\centering
\includegraphics[scale=0.4]{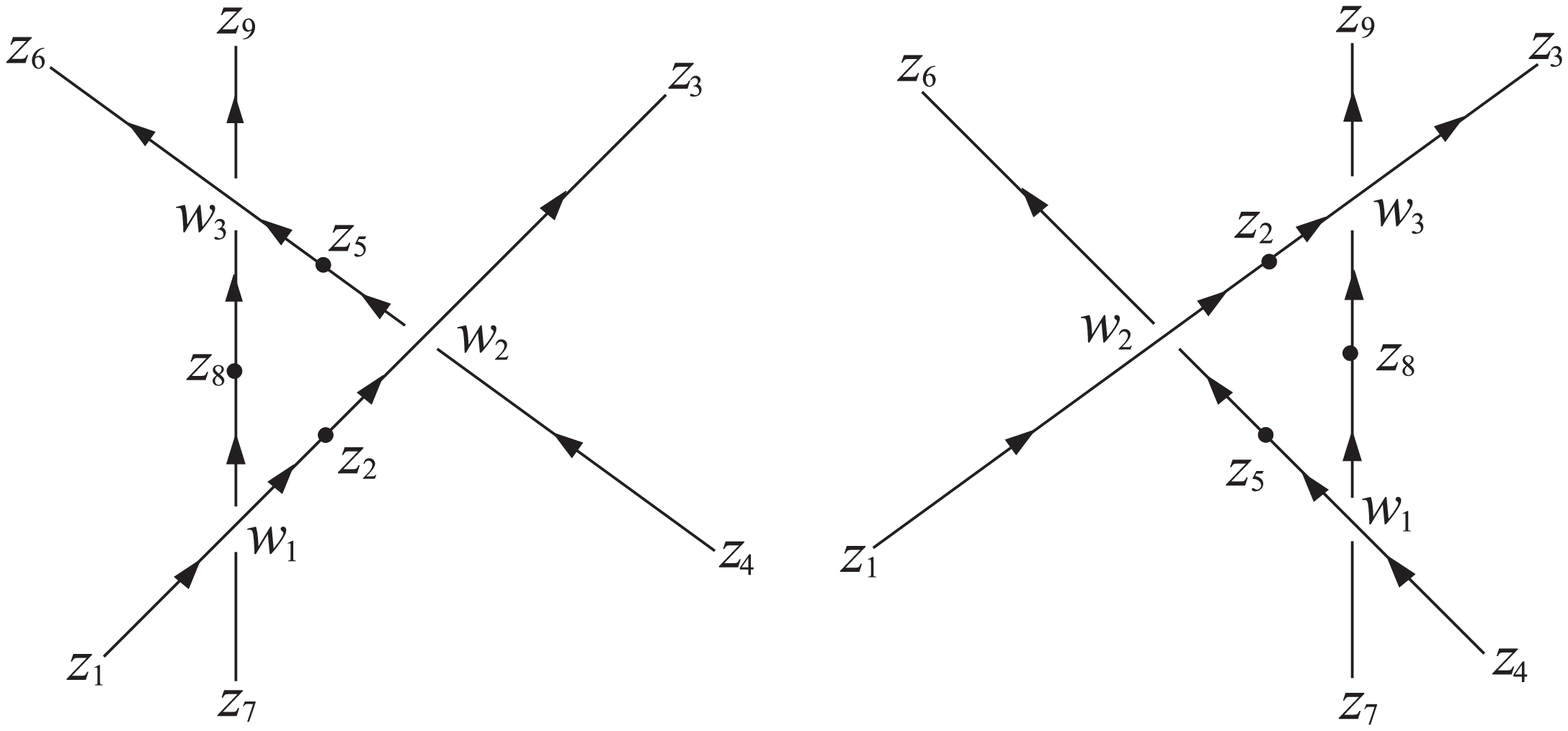}

Fig.5a  \hspace*{5cm} Fig.5b
\end{figure}

On the other hand the definition (\ref{m14}) for the
diagram Fig.5b is given by
\begin{equation}
\begin{array}{rl}
&Tr\langle W(z_7,w_1)W(w_1,z_5)W(z_4,w_1)W(w_1,z_8)\cdot
W(z_8,w_3)W(w_3,z_3)W(z_2,w_3)W(w_3,z_9)\cdot\\
&W(z_5,w_2)W(w_2,z_2)W(z_1,w_2)W(w_2,z_6)
\rangle
\end{array}
\label{m31}
\end{equation}
where the ordering of the three curves is the same as
that in Fig.5a. 
By (\ref{m7a}), (\ref{m8a}), (\ref{m9})  we have that (\ref{m31})
is equal to 
\begin{equation}
Tr\langle W(z_7,z_9)W(z_4,z_5)W(z_1,z_3)R^{-1}W(z_5,z_6)
\rangle
\label{m32}
\end{equation}
Then by reversing the ordering of the curves with end points $z_1$, $z_3$ and with end points $z_4$, $z_6$
respectively we have that both (\ref{m30}) and (\ref{m31})
are equal to
\begin{equation}
Tr\langle W(z_7,z_9)W(z_4,z_6)W(z_1,z_3)R^{-1}
\rangle
\label{m33}
\end{equation}
This shows that Fig.5a is equivalent to Fig.5b and
this gives the Reidemeister move 3.

 From the above examples we see that
knots can be classified by the number $m$ of product of $R$ and $R^{-1}$
matrices. More calculations and examples of the
above knot invariants will be given elsewhere.

\section{Classification of Knots and Links}\label{sec11}

With the knot invariants in the above section we can now
give a classification of knots.
Let $K_1$ and $K_2$ be two knots.
Since the two W-products
$W(z_3,w)W(w,z_2)W(z_1,w)W(w,z_4)$ and $W(z_1,z_2)W(z_3,z_4)$
faithfully represent two oriented pieces of curves which are
crossing or not crossing to each other we have that from the orientation
of a knot $K$ we have that 
the product of
sequence of W-matrices (We may call this product as a generalized
Wilson loop) which are formed according to the orientation
of $K$ in the correlation (\ref{m14}) of defining an invariant of the knot
$K$ faithfully represents the knot $K$.
From this we have that $K_1$ and $K_2$ are topological equivalent if and only if
their generalized Wilson loops  can be transformed to
each other by using (\ref{m7a}), (\ref{m8a}) and (\ref{m9}).
We note that in the above section 
we can derive
the Reidemeister moves by using (\ref{m7a}), (\ref{m8a}) and (\ref{m9}).
Thus this  equivalence of $K_1$ and $K_2$ represented by their 
generalized Wilson loops agrees with the fact that
 $K_1$ and $K_2$ are topologically equivalent if and only
if they can be transformed to each other by the Reidemeister moves.

Then since each knot can be changed to a trivial knot
by applying braiding operation \cite{Mur} 
which is equivalent to (\ref{m7a}), (\ref{m8a}) and (\ref{m9})
it follows that
these new knot invariants can be equivalently transformed to the
form
$Tr R^{-m}\langle W(z_1,z_1)\rangle$
where $m$ is an integer and $Tr \langle W(z_1,z_1)\rangle$ is the
knot invariant for the trivial knot.
This form can also be shown by a direct computation which is similar to
the computation of the invariant of the trivial knot.
We notice that
since this new knot invariant is of the form
$Tr\langle R^{-m}W(z,z)\rangle$ we have that
knots can be completely classified by the power index $m$ of $R$.

Similar to the case of knots we have that the generalized Wilson
loop for a link can faithfully represents this link and that
the correlation (\ref{m14}) of the generalized Wilson loop of this
link is an invariant which can completely
classifies links. 
For the case of link as similar to the case of knot the ordering
of the crossings $W(z_3,w)W(w,z_2)W(z_1,w)W(w,z_4)$ can be given
by following the orientation of each component of a link.
When a component of a link has been traced for one loop 
the crossings 
$W(z_3,w)W(w,z_2)W(z_1,w)W(w,z_4)$ related to this component can be ordered and the pieces $W(w_i,z_j)$ related to these crossings (which may come from other
components of the link) have also been required being ordered.
In the following section we give some examples
to illustrate the formation and computations of this new link invariant.

\section{Examples of New Link Invariants}\label{sec12}

Let us first consider the link in Fig.6a. We may let the two
knots of this link be with $z_1$ and $z_4$ as the initial and final end point respectively. We let the ordering
of these two knots be such that when the $z$ parameter goes
 one loop on one knot then the $z$ parameter for another knot also goes one loop. The correlation
(\ref{m14}) for this link is given by:
\begin{equation}
Tr\langle W(z_3,w_1)W(w_1,z_2)W(z_1,w_1)W(w_1,z_4)\cdot
W(z_4,w_2)W(w_2,z_1)W(z_2,w_2)W(w_2,z_3)
\rangle 
\label{l1a}
\end{equation}
We let the ordering of the $W$-matrices in (\ref{l1a})
be such that $W(z_1,z_2)$ and $W(z_4,z_3)$
start first. Then next $W(z_2,z_1)$ and $W(z_3,z_4)$
follows. Form this ordering we have that (\ref{l1a}) is
equal to:
\begin{equation}
\begin{array}{rl}
&Tr\langle RW(z_1,w_1)W(w_1,z_2)W(z_3,w_1)W(w_1,z_4)\cdot
W(z_4,w_2)W(w_2,z_3)W(z_2,w_2)W(w_2,z_1)R^{-1}
\rangle \\
=&Tr\langle W(z_1,z_2)W(z_3,z_4)
W(z_4,z_3)W(z_2,z_1)
\rangle \\
=&Tr\langle W(z_2,z_2)W(z_3,z_3)
\rangle
\end{array}
\label{l1}
\end{equation}
where we have used (\ref{m7a}) and (\ref{m8a}).
Since by definition (\ref{m14}) we have that $Tr\langle W(z_2,z_2)W(z_3,z_3)
\rangle$ is the knot invariant for two unlinking trivial knots, equation (\ref{l1}) shows that the link in Fig.6a
is topologically equivalent to two unlinking trivial knots.
Similarly we can show that the link in Fig.6b is
topologically equivalent to two unlinking trivial knots.

\begin{figure}[hbt]
\centering
\includegraphics[scale=0.4]{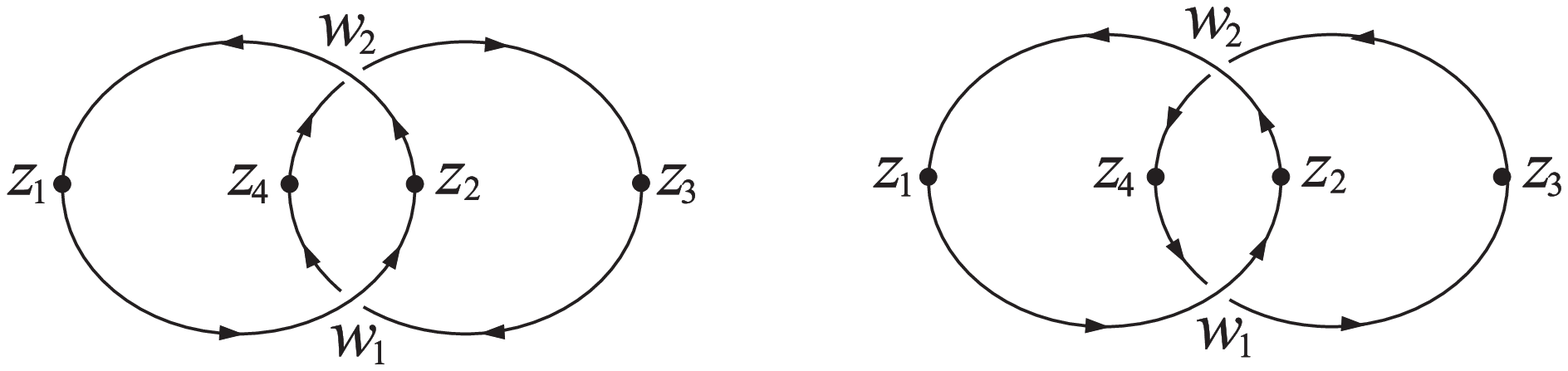}

Fig.6a  \hspace*{3.5cm} Fig.6b
\end{figure}

\begin{figure}[hbt]
\centering
\includegraphics[scale=0.4]{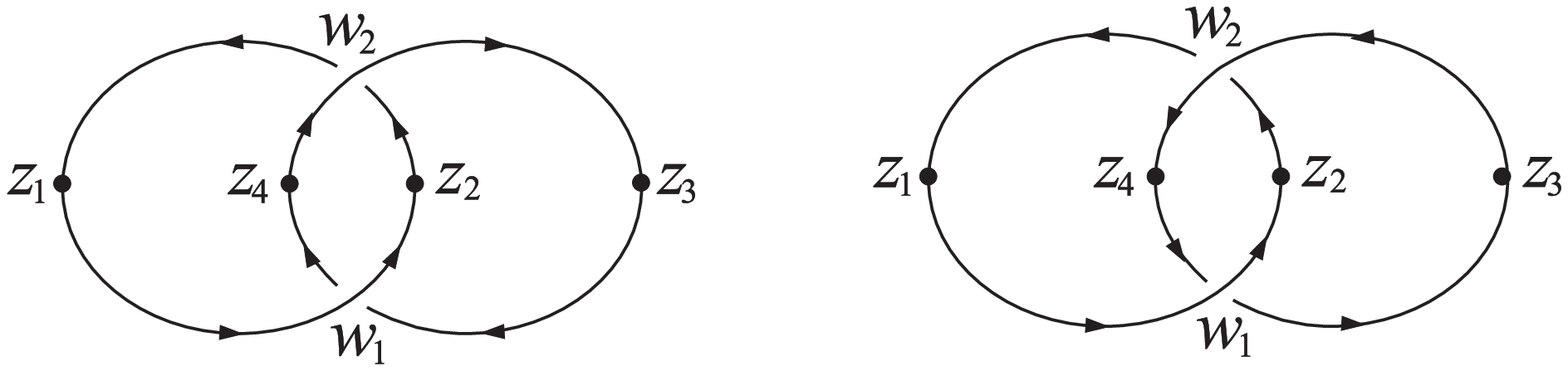}

Fig.7a  \hspace*{3.5cm} Fig.7b
\end{figure}
Let us then consider the Hopf link in Fig.7a. The
correlation (\ref{m14}) for this link is given by:
\begin{equation}
Tr\langle W(z_3,w_1)W(w_1,z_2)W(z_1,w_1)W(w_1,z_4)\cdot
W(z_2,w_2)W(w_2,z_3)W(z_4,w_2)W(w_2,z_1)\rangle
\label{l3}
\end{equation}
The ordering of the $W$-matrices in (\ref{l3})
is such that $W(z_1,z_2)$ starts first and $W(z_3,z_4)$
follows it. Then next  we let $W(z_2,z_1)$ starts first
and $W(z_4,z_3)$ follows it. The ordering is such that
when the $z$ parameter goes
one loop in one knot of the link we have that the
$z$ parameter also goes one loop on the other knot.
From the ordering we have that (\ref{l3}) is equal to:
\begin{equation}
\begin{array}{rl}
&Tr\langle RW(z_1,w_1)W(w_1,z_2)W(z_3,w_1)W(w_1,z_4)\cdot\\
&W(z_2,w_2)W(w_2,z_1)W(z_4,w_2)W(w_2,z_3)R^{-1}\rangle \\
=&Tr\langle W(z_1,z_2)W(z_3,z_4)
W(z_2,z_1)W(z_4,z_3)\rangle
\end{array}
\label{l4}
\end{equation}
Then let us consider the following correlation:
\begin{equation}
Tr\langle 
R^{-2}W(z_3,w_1)W(w_1,z_2)W(z_1,w_1)W(w_1,z_4)\cdot
W(z_2,w_2)W(w_2,z_3)W(z_4,w_2)W(w_2,z_1)\rangle 
\label{l5}
\end{equation}
We let the ordering of the $W$-matrices in (\ref{l5})
be such that $W(z_1,z_2)$ starts first and $W(z_4,z_3)$
follows it. Then next $W(z_2,z_1)$ starts first and
$W(z_3,z_4)$ follows it. From the ordering we have that
(\ref{l5}) is equal to:
\begin{equation}
\begin{array}{rl}
&Tr\langle 
R^{-2}RW(z_1,w_1)W(w_1,z_2)W(z_3,w_1)W(w_1,z_4)\cdot
W(z_2,w_2)W(w_2,z_1)W(z_4,w_2)W(w_2,z_3)R\rangle \\
=&Tr\langle W(z_1,z_2)W(z_3,z_4)
W(z_2,z_1)W(z_4,z_3)\rangle
\end{array}
\label{l6}
\end{equation}
On the other hand from the ordering of (\ref{l5}) we
have that (\ref{l5}) is equal to:
\begin{equation}
\begin{array}{rl}
&Tr\langle R^{-2}W(z_3,w_1)RW(z_1,w_1)W(w_1,z_2)R^{-1}\\
&W(w_1,z_4)W(z_2,w_2)
RW(z_4,w_2)W(w_2,z_3)R^{-1}W(w_2,z_1)\rangle \\
=&Tr\langle R^{-2}W(z_3,w_1)RW(z_1,z_2)
R^{-1}W(w_1,z_4)W(z_2,w_2)R
W(z_4,z_3)R^{-1}W(w_2,z_1)\rangle \\
=&Tr\langle R^{-2}W(z_3,w_1)RW(z_1,z_2)
W(z_2,w_2)W(w_1,z_4)W(z_4,z_3)R^{-1}W(w_2,z_1)\rangle \\
=&Tr\langle R^{-2}W(z_3,w_1)
RW(z_1,w_2)W(w_1,z_3)R^{-1}W(w_2,z_1)\rangle \\
=&Tr\langle R^{-2}W(z_3,w_1)
W(w_1,z_3)W(z_1,w_2)W(w_2,z_1)\rangle \\
=&Tr\langle R^{-2}W(z_3,z_3)W(z_1,z_1)\rangle 
\end{array}
\label{l7}
\end{equation}
where we have repeatly used (\ref{m9}).
From (\ref{l4}), (\ref{l6}) and (\ref{l7}) we have that
the knot invariant for the Hopf link in Fig.7a is given by:
\begin{equation}
Tr\langle R^{-2}W(z_3,z_3)W(z_1,z_1)\rangle
\label{l8}
\end{equation}

Then let us consider the Hopf link in Fig.7b. The correlation for this link is given by
\begin{equation}
Tr\langle W(z_4,w_1)W(w_1,z_2)W(z_1,w_1)W(w_1,z_3)\cdot
W(z_2,w_2)W(w_2,z_4)W(z_3,w_2)W(w_2,z_1)\rangle
\label{l9}
\end{equation}
By a derivation which is dual to the above derivation
for the Hopf link in Fig.7a we have that (\ref{l9})
is equal to
\begin{equation}
Tr\langle R^{2}W(z_4,z_4)W(z_1,z_1)\rangle
\label{l10}
\end{equation}
We see that the invariants for the above two Hopf links
are different. This agrees with the fact that these
two links are not topologically equivalent.

We can extend the above computations to other links.
As examples let us consider the linking of two trivial
knot with linking number $n$. Similar to the above
computations we have that this link which analogous to
the Hopf link in Fig.7a is with an invariant equals
to $Tr\langle R^{-2n}W(z_4,z_4)W(z_1,z_1)\rangle$. Also
for this link which analogous to to the Hopf link
in Fig.7b is with an invariant equals to
$Tr\langle R^{2n}W(z_4,z_4)W(z_1,z_1)\rangle$.

More calculations and examples
 of these new link invariants will be given elsewhere.

\section{A Classification Table of Knots}\label{sec13}

In this section
let us give some arguments to
determine the new knot invariants of prime knots
and nonprime knots. We have shown that this new 
invariant of each knot is of the form
$Tr\langle R^{-m}W(z_1,z_1)\rangle$
where $m$ is an integer. This power index $m$ can be
regarded as a measure of the complexity of a knot. Let us determine $m$ for prime and nonprime knots. We need only to
determine $m$ for knots with positive $m$ since the
corresponding mirror image will have negative $m$
if the mirror image is not equivalent to the corresponding
knot.
We have shown that the invariants for links of two
trivial knots with linking number $n$ as in Fig.7a
are with the
product of $R$ of the form $R^{-2n}$ with $m=2n$.
This $m$ is an even number. Now if we insert these links
into the knot table of prime knots we see that the
prime knots must with $m$ being an odd number (We may refer
to the knot table in \cite{Ka}. We may have nonprime knots
which are not in this knot table having the same $m$ as
these links. In this case our new link invariants still
can distinguish them because the corresponding link
invariants are of two-loop form while the corresponding
knot invariants are of one-loop form).
 
Then we expect that for prime knots $m$ is an odd
prime number. Computations show that for the knot $\bf 3_1$ we
have $m=1$, for the knot $\bf 4_1$ we have $m=3$. Then from some arguments on the effect of $R$ we
should have that for the knot $\bf 5_1$ we have $m=5$
and for the knot $\bf 5_2$ we have $m=7$.
Then how about the knot $\bf 6_1$? We have that the numbers
from $1$ to $8$ are ocupied by knots and links of two trivial knots. Then $10$ is ocupied by the link of two
trivial knots with linking number $5$. Thus for the knot
$\bf 6_1$ we should have $m=9$ or $m=11$. From some
argument on the effect of $R$ we should have $m=11$ for
$\bf 6_1$.
Then is there a knot with $m=9$?

Let us first consider the granny knot (or the square knot) which is a nonprime knot composed with the knot $\bf 3_1$ and its mirror image. This square knot
has $6$ crossings and $4$ alternating crossings and thus its complexity
which is measured by the power index $m$ of $R$ is
less than that of $\bf 5_1$ which is with $5$ alternating crossings. Thus this granny knot is with $m=4$ ($m=3$
has been occupied by $\bf 4_1$).
Let us denote this granny knot by ${\bf 3_1\star 3_1}$
where $\star$ denotes the connected sum of two knots such that 
the resulting total number of alternating crossings
is equal to the total number of alternating crossings of the two knots minus $2$. 

Then let us consider the
reef knot which is a nonprime knot composed with two identical knots $\bf 3_1$. This knot has $6$ alternating crossings which is equal to the total number of crossings
as that of $\bf 6_1$. Since this knot is nonprime its
complexity is less than that of $\bf 6_1$ where the
complexity may be measured by the power index $m$ of $R$.
Thus this reef knot is with power index $m$ less than
$11$.
Then if we also regard the total number of alternating
crossings of a knot as a way to measure the complexity
of a knot we have that the power index $m$ of this
granny knot is greater than $5$ since $\bf 5_1$ is
with $5$ alternating crossings and with $m=5$.
Let us denote this reef knot by ${\bf 3_1\times 3_1}$
where $\times$ denotes the connected sum for two knots
such that the resulting total number of alternating
crossings is equal to the total number of alternating
crossings of the two knots.

Now let us look for knots with $5$ or $6$ alternating crossings.
Let us consider the nonprime knot 
$\bf 3_1\star 4_1$ composed with
a knot $\bf 3_1$ and a knot $\bf 4_1$ with $7$ crossings
and $5$ alternating crossings. In this case we have that
the power index $m$ of $\bf 3_1\star 4_1$ should be greater than that
of  $\bf 5_1$ which is exactly with  $5$ alternating crossings
since $\bf 3_1\star 4_1$ in addtion has $7$ crossings.
Then the power index $m$ of $\bf 3_1\star 4_1$ 
should be less than that of $\bf 5_2$ which is also with  $5$ alternating crossings but these crossings are arranged
in a more complicated way which is an effect of $R^2$
such that $\bf 5_2$ is with $m=7$. Thus 
$\bf 3_1\star 4_1$ is with $m=6$.

Then we consider the nonprime knot 
$\bf 3_1\star (3_1\star 3_1)$ with $9$
crossings and $5$ alternating crossings. The power index
of this knot should be 
less than that of $\bf 3_1\times 3_1$ since 
it is with $6$ alternating crossings.
Then the power index
of $\bf 3_1\times 3_1$  should be 
greater than that of $\bf 5_2$ since it has in addition
$9$ crossings which would be enough for
a greater power index $m$. Then we have that
$\bf 3_1\star (3_1\star 3_1)$ is with $m=8$ since
it has $9$ crossings and $5$ alternating crossings
and thus is with the same complexity as $\bf 5_2$.
Then since $\bf 5_2$ can not have $m=8$ we thus have
that $\bf 3_1\star (3_1\star 3_1)$ is with $m=8$.

Then we consider the nonprime knot 
$\bf 3_1\star 5_1$ with $8$ crossings and $6$
alternating crossings. The power index $m$ of this
knot should be greater than that of $\bf 3_1\times 3_1$
which has exactly $6$ alternating crossings.
Then the power index $m$ of this
knot should be less than that of $\bf 6_1$ which also has
$6$ alternating crossings but these crossings are arranged
in a more complicated way with an effect of $R^4$ from
that of $5_2$. Thus $\bf 3_1\star 5_1$ is with $m=10$
and finally we have that $\bf 3_1\times 3_1$ is with
power index $m=9$.

Thus for $m$ from $1$ to $11$ we have fill in a suitable
knot with power index $m$ (except the case $m=2$ which
is filled in with the Hopf link) such that odd prime numbers
are filled with prime knots.
In a similar way we may determine the power index $m$
of other prime and nonprime knots. We list the results
up to $m=2^5$ in a form of table.

\begin{displaymath}
\begin{array}{|c|c|c|c|} \hline
\mbox{Type of Knot}& \mbox{ Power \, Index}\: m 
 &\mbox{Type of Knot}& \mbox{ Power \, Index}\: m 
\\ \hline
{\bf 3_1} & 1 & {\bf 6_3} &  17\\ \hline

\mbox{Hopf link} &  2 
&  {\bf 3_1\times 4_1} &  18 \\ \hline

{\bf 4_1} &  3 & {\bf 7_1} &  19 \\ \hline

{\bf 3_1\star 3_1} &  4 & 
{\bf 4_1\star 5_1} &  20
\\ \hline

{\bf 5_1} & 5 & {\bf 4_1\star(3_1\star 4_1) } &  21
\\ \hline

{\bf 3_1\star 4_1} & 6 & {\bf 4_1\star 5_2} & 22 \\ \hline

{\bf 5_2} &  7 & {\bf 7_2} & 23 \\ \hline

{\bf 3_1\star 3_1\star 3_1} &  8 & 
{\bf 3_1\star (3_1\times 3_1)}& 24 \\ \hline

{\bf 3_1\times 3_1} & 9 &
{\bf 3_1\star (3_1\star  5_1)}& 25 \\ \hline

{\bf 3_1\star 5_1} &  10 & {\bf 3_1\star 6_1}
& 26 \\ \hline

{\bf 6_1} &  11 & {\bf 3_1\star (3_1\star  5_2)}& 27
\\ \hline

{\bf 3_1\star 5_2} &  12 & {\bf 3_1\star 6_2} & 28 
\\ \hline

{\bf 6_2} &  13 & {\bf 7_3} & 29 \\ \hline

{\bf 4_1\star 4_1} &  14 & 
{ \bf 3_1\star (3_1\star 3_1)\star 4_1} & 30 \\ \hline

{\bf 4_1\star (3_1\star 3_1)} & 15 & {\bf 7_4} & 31 
\\ \hline

{ \bf (3_1\star 3_1)\star (3_1\star 3_1)} & 16 &
{\bf 3_1\star (3_1\star 3_1)\star (3_1\star 3_1)} &  32 
\\ \hline
\end{array}
\end{displaymath}

From this table we see that comparable nonprime knots (in a sense from the table) are grouped
in each of the intervals between two prime numbers. It is interesting that in 
each interval nonprime numbers are one-to-one  filled with the comparable nonprime knots
while prime numbers are filled with prime knots. This grouping property of classification reflects that the new knot invariants
(and hence the power index $m$) give a classification of
knots.

Let us find out some rules for the whole classification
table.
We have shown that
even numbers can not be filled with prime knots.
Thus even numbers (except $2$) can only be filled with nonprime knots.
On the other hand each odd nonprime number is between
two even numbers which are power indexes of nonprime
knots. Thus the knot corresponding to this odd number
is in the same group of these two nonprime knots and
is comparable with these two nonprime knots and thus
must also be a nonprime knot.
From this we have that nonprime numbers are power indexes of nonprime knot. Similarly by this grouping property
we have that odd prime numbers are power indexes of
prime knots.  

It is interesting to note that from the above knot table
We see that in the $\star$ product 
the knot ${\bf 3_1}$ plays the role of the number $2$ 
in the usual multiplication of numbers. Thus the $\star$ product (or the connected sum) is a kind of  multiplication corresponding to the usual multiplication
of numbers. However the general rule for this multiplication
is rather complicated. This reflects the fact that
the numbers $m$ are the power indexes of $R$ which are with simple rule for
the addition ( and not for multipication).
From this generalized multiplication (which is the connected sum) of knots which corresponds to the usual
multiplication of the power indexes $m$ we also have that
prime knots are with odd prime numbers as power indexes
and nonprime knots are with nonprime numbers as power indexes.

Thus we have a classification table of knots such that
each prime knot corresponds to a prime power index $m$
and each nonprime knot corresponds to a nonprime power index $m$.
More computations to verify this knot table shall be given alsewhere.

\section{Conclusion}

In this paper from a new quantum field model we have derived a conformal field theory and a quantum group
structure for generalized Wilson loops from which we can derive
the Jones polynomial and new knot and link invariants which extend
the Jones polynomial. We show that these new invariants can completely classify knots and links. These
new invariants are in terms of the monodromy $R$
of two
Knizhnik-Zamolochikov equations which are dual to each other. In the case of knots these new invariants can be written in the form
$Tr\langle R^{-m}W(z,z)\rangle$ from which
we may classify
knots with the power index $m$ of $R$. A classification table of knots can then be formed where prime knots
are classified with odd prime numbers $m$ and nonprime
knots are classified with nonprime numbers $m$.

\end{document}